\documentclass{article}

\usepackage{amsmath,amssymb,amsthm}

\newtheorem{theorem}{Theorem}
\newtheorem{lemma}[theorem]{Lemma}
\newtheorem{conjecture}[theorem]{Conjecture}

\begin{document}

\title{Smarandache Sequences: Explorations
       and Discoveries with a Computer Algebra System\footnote{To appear
       in the 2003 December issue of the \emph{Smarandache Notions Journal}.}}

\author{Paulo~D.~F.~Gouveia \\
        \texttt{pgouveia@ipb.pt} \\
        Technology and Manag. School\\
        Bragan\c{c}a Polytechnic Institute\\
        5301-854 Bragan\c{c}a, Portugal
        \and
        Delfim~F.~M.~Torres \\
        \texttt{delfim@mat.ua.pt} \\
        Department of Mathematics \\
        University of Aveiro \\
        3810-193 Aveiro, Portugal}

\date{}

\maketitle


\begin{abstract}
We study Smarandache sequences
of numbers, and related problems,
via a Computer Algebra System.
Solutions are discovered, and some
conjectures presented.
\end{abstract}


\bigskip

\noindent \textbf{Mathematics Subject Classification 2000.} 11B83, 11-04, 68W30.

\bigskip


\noindent \textbf{Keywords.} Smarandache square/cube/prime digital subsequence,
Smarandache Lucas/Fibonacci partial digital subsequence, Smarandache
odd/even/
prime/Fibonacci sequence, Smarandache function, Smarandache
2-2-additive/2-2-subtractive/3-3-additive relationship, Smarandache partial
perfect additive sequence, Smarandache bad numbers, Smarandache prime
conjecture, mathematical experimentation.


\section{Introduction}

After a good look on the \emph{Mathematics Unlimited---2001 and Beyond}
\cite{MR2002d:00004a}, which addresses the question of the future of Mathematics
in the new millennium, it is impossible not to get the deep impression
that \emph{Computing} will be an integral part of many branches
of Mathematics. If it is true that in the {XX}st century Mathematics
has contributed, in a fundamental way, to technology, now, in the
{XXI}st century, the converse seems to be also a possibility. For perspectives
on the role of Computing in Mathematics (and the other way around)
see \cite{MR1852153,MR2000m:68186,MR1754789}.

Many powerful and versatile Computer Algebra Systems
are available nowadays, putting at our disposal
sophisticated environments of mathematical and
scientific computing.
They comprise both numerical and symbolic computation
through high-level and expressive languages,
close to the mathematical one. A large quantity of mathematical
knowledge is already available in these scientific systems,
providing efficient mathematical methods to perform the
desired calculations. This has two important implications: they spare
one a protracted process of programming and debugging,
so common to the more conventional computer languages;
they permit us to write few lines of code, and simpler programs,
more declarative in nature. Our claim is that explorations
with such tools can develop intuition, insight, and better qualitative
understanding of the nature of the problems. This can greatly assist
the proof of mathematical results (see an example in Section
\S~\ref{SubSecNR} below).

It is our aim to show that computer-assisted algebra
can provide insight and clues to some open questions related
to special sequences in Number Theory. Number Theory has
the advantage of being easily amenable to computation and experimentation.
Explorations with a Computer Algebra System will allow us to
produce results and to formulate conjectures.
We illustrate our approach with the mathematics \textsf{Maple} system
(all the computational processing was carried with \textsf{Maple} version 8,
on an AMD Athlon(TM) 1.66 GHz machine), and with some of the problems proposed
by the Romanian mathematician Florentin Smarandache.

\textsf{Maple} was originated more than two decades ago, as a
project of the Symbolic Computation Group of the University of
Waterloo, Ontario. It is now a registered trademark product of
Waterloo Maple Inc. We refer the reader to
\cite{MR1823705,LivroRosenEDMM} for a gentle introduction to
\textsf{Maple}. For a good account on the Smarandache collection
of problems, and for a biography of F.~Smarandache, see
\cite{MR97m:01073}.

We invite and exhort readers to convert our
mathematical explorations in the language of their favorite
Computer Algebra System; to optimize the algorithms
(we have followed the didactic approach, without any attempt
of code optimization); and to obtain the results for themselves.
The source be with you.


\section{Smarandache Digital Subsequences}

We begin by considering sequences of natural numbers
satisfying some given property together
with all their digits.

\subsection{Smarandache p-digital subsequences}
\label{SubSecNR}

We are interested in the following
\emph{Smarandache p-digital subsequences}.
Let $p \ge 2$. From the sequence $\left\{n^p\right\}$,
$n \in \mathbb{N}_0$, we select those terms whose
digits are all perfect $p$-powers. For $p = 2$ we
obtain the \emph{Smarandache square-digital subsequence}:
we select only those terms of the sequence
$\left\{n^2\right\}_{n=0}^{\infty}$
whose digits belong to the set $\left\{0,1,4,9\right\}$.
With the \textsf{Maple} definitions
\begin{verbatim}
> pow := (n,p) -> seq(i^p,i=0..n):
> perfectPow := (n,p) -> evalb(n = iroot(n,p)^p):
> digit := (n,num) -> irem(iquo(num,10^(length(num)-n)),10):
> digits := n -> map(digit,[$1..length(n)],n):
> digPerfectPow :=
>   (n,p) -> evalb(select(perfectPow,digits(n),p) = digits(n)):
\end{verbatim}
the Smarandache square-digital subsequence is easily obtained:
\begin{verbatim}
> ssds := n -> select(digPerfectPow,[pow(n,2)],2):
\end{verbatim}
We now ask for all the terms of
the Smarandache square-digital subsequence
which are less or equal than $10000^2$:
\begin{verbatim}
> ssds(10000);
\end{verbatim}
\begin{gather*}
[0, 1, 4, 9, 49, 100, 144, 400, 441, 900, 1444, 4900, 9409, 10000, 10404,
11449,\\
14400, 19044, 40000, 40401, 44100, 44944, 90000, 144400, 419904, 490000,\\
491401, 904401, 940900, 994009, 1000000, 1004004, 1014049, 1040400,\\
1100401, 1144900, 1440000, 1904400, 1940449, 4000000, 4004001, 4040100,\\
4410000, 4494400, 9000000, 9909904, 9941409, 11909401, 14010049, 14040009,\\
14440000, 19909444, 40411449, 41990400, 49000000, 49014001, 49140100,\\
49999041, 90440100, 94090000, 94109401, 99400900, 99940009, 100000000]
\end{gather*}

In \cite{MR1764401,MR1764402} one finds the following question:
\begin{quotation}
``Disregarding the square numbers of the form $N \times 10^{2k}$,
$k \in \mathbb{N}$, $N$ also a perfect square number,
how many other numbers belong to the Smarandache square-digital subsequence?''
\end{quotation}
From the obtained 64 numbers of the Smarandache square-digital subsequence,
one can see some interesting patterns from which one easily guess
the answer.
\begin{theorem}
\label{th:mr}
There exist an infinite number of terms on the Smarandache square-digital
subsequence which are not of the form $N \times 10^{2k}$,
$k \in \mathbb{N}$, $N$ a perfect square number.
\end{theorem}

Theorem~\ref{th:mr} is a straightforward consequence
of the following Lemma.
\begin{lemma}
\label{lm:mr}
Any number of the form $\left(10^{k+1}+4\right)\times 10^{k+1} + 4$,
$k \in \mathbb{N}_0$ ($144$, $10404$, $1004004$, $100040004$, ...),
belong to the Smarandache square-digital subsequence.
\end{lemma}

\begin{proof}
Lemma~\ref{lm:mr} follows by direct calculation:
\begin{equation*}
\left(10^{k+1}+2\right)^2 = \left(10^{k+1}+4\right)\times 10^{k+1} + 4 \, .
\end{equation*}
\end{proof}

We remark that from the analysis of the list
of the first 64 terms of the Smarandache
square-digital subsequence, one easily finds
other possibilities to prove Theorem~\ref{th:mr},
using different but similar assertions than the one
in Lemma~\ref{lm:mr}. For example, any number of the form
$\left(10^{k+2}+14\right) \times 10^{k+2} + 49$,
$k \in \mathbb{N}_0$ ($11449$, $1014049$, $100140049$, ...),
belong to the Smarandache square-digital subsequence:
\begin{equation*}
\left(10^{k+2}+7\right)^2 = \left(10^{k+2}+14\right) \times 10^{k+2} + 49 \, .
\end{equation*}
Other possibility, first discovered in \cite{MaohuaRSSDS}, is to use the
pattern $\left(4 \times 10^{k+1}+4\right) \times 10^{k+1} + 1$,
$k \in \mathbb{N}_0$ ($441$, $40401$, $4004001$, ...),
which is the square of $2 \times 10^{k+1} + 1$.

Choosing $p=3$ we obtain
the \emph{Smarandache cube-digital subsequence}.
\begin{verbatim}
> scds := n -> select(digPerfectPow,[pow(n,3)],3):
\end{verbatim}
Looking for all terms of the Smarandache cube-digital subsequence
which are less or equal than $10000^3$
we only find the trivial ones:
\begin{verbatim}
> scds(10000);
\end{verbatim}
\begin{equation*}
[0, 1, 8, 1000, 8000, 1000000, 8000000, 1000000000, 8000000000, 1000000000000]
\end{equation*}
We offer the following conjecture:
\begin{conjecture}
All terms of the Smarandache cube-digital subsequence
are of the form $D \times 10^{3k}$ where
$D \in \left\{0,1,8\right\}$ and $k \in \mathbb{N}_0$.
\end{conjecture}

Many more Smarandache digital subsequences have been introduced
in the literature. One good example is the Smarandache
prime digital subsequence, defined as the sequence of prime numbers
whose digits are all primes (see \cite{MR1764402}).

Terms of the Smarandache prime digital subsequence are easily
discovered with the help of the \textsf{Maple} system. Defining
\begin{verbatim}
> primeDig := n -> evalb(select(isprime,digits(n)) = digits(n)):
> spds := n -> select(primeDig,[seq(ithprime(i),i=1..n)]):
\end{verbatim}
we find that 189 of the first 10000 prime numbers
belong to the Smarandache prime digital subsequence:
\begin{verbatim}
> nops(spds(10000));
\end{verbatim}
\begin{equation*}
189
\end{equation*}


\subsection{Smarandache p-partial digital subsequences}

The \emph{Smarandache p-partial digital subsequence} is
defined by scrolling through a given sequence $\{a_n\}$,
$n \geq 0$, defined by some property $p$, and selecting
only those terms which can be partitioned in groups
of digits satisfying the same property $p$ (see \cite{MR1764401}).
For example, let us consider $\{a_n\}$ defined by the
recurrence relation $a_n = a_{n-1} + a_{n-2}$.
One gets the \emph{Lucas sequence} by choosing the initial
conditions $a_0 = 2$ and $a_1 = 1$; the \emph{Fibonacci sequence}
by choosing $a_0 = 0$ and $a_1 = 1$.
The Smarandache \emph{Lucas}-partial digital subsequence
and the Smarandache \emph{Fibonacci}-partial digital subsequence
are then obtained selecting from the respective sequences
only those terms $n$ for which there exist a partition of the digits
in three groups ($n = g_1 g_2 g_3$) with the sum of the first two groups
equal to the third one ($g_1 + g_2 = g_3$).

In \cite{MR1764401,oai:arXiv.org:math/0010132,PropNumb}
the following questions are formulated:
\begin{quotation}
``Is 123 ($1+2=3$) the only Lucas number that verifies
a Smarandache type partition?''
\end{quotation}
\begin{quotation}
``We were not able to find any Fibonacci number verifying a Smarandache
type partition, but we could not investigate large numbers; can you?''
\end{quotation}

Using the following procedure, we can verify if a certain number $n$ fulfills
the necessary condition to belong to the Smarandache
\emph{Lucas}/\emph{Fibonacci}-partial digital subsequence, \textrm{i.e.},
if $n$ can be divided in three digit groups, g1g2g3, with g1+g2=g3.

\begin{verbatim}
> spds:=proc(n)
>   local nd1, nd2, nd3, nd, g1, g2, g3:
>   nd:=length(n);
>   for nd3 to nd-2 do
>     g3:=irem(n,10^nd3);
>     if length(g3)*2>nd then break; fi;
>     for nd1 from min(nd3,nd-nd3-1) by -1 to 1 do
>       nd2:=nd-nd3-nd1;
>       g1:=iquo(n,10^(nd2+nd3));
>       g2:=irem(iquo(n,10^nd3), 10^nd2);
>       if g2>=g3 then break;fi;
>       if g1+g2=g3 then printf("%d (%d+%d=%d)\n",n,g1,g2,g3);fi;
>     od;
>   od:
> end proc:
\end{verbatim}

Now, we can compute the first $n$ terms of the Lucas sequence, using the
procedure below.

\begin{verbatim}
> lucas:=proc(n)
>   local L, i:
>   L:=[2, 1]:
>   for i from 1 to n-2 do L:=[L[],L[i]+L[i+1]]:od:
> end proc:
\end{verbatim}
With $n=20$ we get the first twenty Lucas numbers
\begin{verbatim}
> lucas(20);
\end{verbatim}
\begin{equation*}
[2, 1, 3, 4, 7, 11, 18, 29, 47, 76, 123, 199, 322, 521, 843, 1364, 2207, 3571,
5778, 9349]
\end{equation*}
Let $L$ be the list of the first $6000$ terms of the Lucas sequence:
\begin{verbatim}
> L:=lucas(6000):
\end{verbatim}
\verb|(elapsed time: 1.9 seconds)|
\footnote{The most significant time calculations are showed, in
order to give an idea about the involved computation effort.}
\\

It is interesting to remark that the $6000^{th}$ element has 1254 digits:

\begin{verbatim}
> length(L[6000]);
\end{verbatim}
\begin{equation*}
1254
\end{equation*}

The following \textsf{Maple} command permit us to check which
of the first $3000$ elements belong to a Smarandache
\emph{Lucas}-partial digital subsequence.

\begin{verbatim}
> map(spds, L[1..3000]):
\end{verbatim}
\begin{verbatim}
123 (1+2=3)
20633239 (206+33=239)
\end{verbatim}
\begin{verbatim}
(elapsed time: 7h50m)
\end{verbatim}

As reported in \cite{fibonacciSubSeq}, only two of the first $3000$
elements of the Lucas sequence verify a Smarandache type partition: the
$11^{th}$ and $36^{th}$ elements.

\begin{verbatim}
> L[11], L[36];
\end{verbatim}
\begin{equation*}
123, 20633239
\end{equation*}

We now address the following question:
Which of the next $3000$ elements of the Lucas sequence belong to a
Smarandache \emph{Lucas}-partial digital subsequence?

\begin{verbatim}
> map(spds, L[3001..6000]):
\end{verbatim}
\begin{verbatim}
(elapsed time: 67h59m)
\end{verbatim}

The answer turns out to be \emph{none}:
no number, verifying a Smarandache type partition, was found between the
$3001^{th}$ and the $6000^{th}$ term of the Lucas sequence.
\\
\\
The same kind of analysis is easily done for the Fibonacci sequence.
We compute the terms of the Fibonacci sequence using the
pre-defined function \texttt{fibonacci}:

\begin{verbatim}
> with(combinat, fibonacci):
> [seq(fibonacci(i), i=1..20)];
\end{verbatim}
\begin{equation*}
[1, 1, 2, 3, 5, 8, 13, 21, 34, 55, 89, 144, 233, 377, 610, 987, 1597, 2584,
4181, 6765]
\end{equation*}

Although the $6000^{th}$ Fibonacci number is different
from the $6000^{th}$ Lucas number
\begin{verbatim}
> evalb(fibonacci(6000) = L[6000]);
\end{verbatim}
\begin{equation*}
false
\end{equation*}
they have the same number of digits
\begin{verbatim}
> length(fibonacci(6000));
\end{verbatim}
\begin{equation*}
1254
\end{equation*}

In order to identify which of the first $3000$ Fibonacci numbers
belong to the Smarandache
\emph{Fibonacci}-partial digital subsequence, we execute the following
short piece of \textsf{Maple} code:

\begin{verbatim}
> map(spds, [seq(fibonacci(i), i=1..3000)]):
\end{verbatim}
\begin{verbatim}
832040 (8+32=040)
\end{verbatim}
\begin{verbatim}
(elapsed time: 8h32m)
\end{verbatim}

This is in consonance with the result
reported in \cite{fibonacciSubSeq}: only one number, among the first $3000$
numbers of the Fibonacci sequence, verifies a Smarandache type partition
-- the $30^{th}$ one.

\begin{verbatim}
> fibonacci(30);
\end{verbatim}
\begin{equation*}
832040
\end{equation*}

As before, with respect to the Lucas sequence,
we now want to know which of the next $3000$ numbers
of the Fibonacci sequence belong to the Smarandache
\emph{Fibonacci}-partial digital subsequence.

\begin{verbatim}
> map(spds, [seq(fibonacci(i), i=3001..6000)]):
\end{verbatim}
\begin{verbatim}
(elapsed time: 39h57m)
\end{verbatim}

Similarly to the Lucas case, no number, verifying a Smarandache
type partition, was found between the $3001^{th}$ and the
$6000^{th}$ term of the Fibonacci sequence.


\section{Smarandache Concatenation-Type Sequences}

Let $\left\{a_n\right\}$, $n \in \mathbb{N}$,
be a given sequence of numbers. The Smarandache
concatenation sequence associated to $\left\{a_n\right\}$
is a new sequence $\left\{s_n\right\}$
where $s_n$ is given by the concatenation of all the terms
$a_1$, $\ldots$, $a_n$. The concatenation operation between
two numbers $a$ and $b$ is defined as follows:
\begin{verbatim}
> conc := (a,b) -> a*10^length(b)+b:
\end{verbatim}
In this section we consider four different Smarandache
concatenation-type subsequences: the odd, the even, the prime,
and the Fibonacci one.
\begin{verbatim}
> oddSeq   := n -> select(type,[seq(i,i=1..n)],odd):
> evenSeq  := n -> select(type,[seq(i,i=1..n)],even):
> primeSeq := n -> [seq(ithprime(i),i=1..n)]:
> with(combinat, fibonacci):
> fibSeq   := n -> [seq(fibonacci(i),i=1..n)]:
> # ss = Smarandache Sequence
> ss := proc(F,n)
>   local L, R, i:
>   L := F(n):
>   R := array(1..nops(L)): R[1] := L[1]:
>   for i from 2 while i <= nops(L) do
>     R[i]:=conc(R[i-1],L[i]):
>   end do:
>   evalm(R):
> end proc:
\end{verbatim}
Just to illustrate the above definitions,
we compute the first five terms of the Smarandache
odd, even, prime, and Fibonacci sequences:
\begin{verbatim}
> ss(oddSeq,10);
\end{verbatim}
\begin{equation*}
[1, 13, 135, 1357, 13579]
\end{equation*}
\begin{verbatim}
> ss(evenSeq,10);
\end{verbatim}
\begin{equation*}
[2, 24, 246, 2468, 246810]
\end{equation*}
\begin{verbatim}
> ss(primeSeq,5);
\end{verbatim}
\begin{equation*}
[2, 23, 235, 2357, 235711]
\end{equation*}
\begin{verbatim}
> ss(fibSeq,5);
\end{verbatim}
\begin{equation*}
[1, 11, 112, 1123, 11235]
\end{equation*}
Many interesting questions appear when one try to find
numbers among the terms of a Smarandache concatenation-type
sequence with some given property. For example,
it remains an open question to understand how many primes
are there in the odd, prime, or Fibonacci sequences.
Are they infinitely or finitely in number?
The following procedure permit us to find
prime numbers in a certain Smarandache sequence.
\begin{verbatim}
> ssPrimes := proc(F,n)
>   local ar, i:
>   ar := select(isprime,ss(F,n)):
>   convert(ar,list):
> end proc:
\end{verbatim}
There are five prime numbers in the first fifty
terms of the Smarandache odd sequence;
\begin{verbatim}
> nops(ssPrimes(oddSeq,100));
\end{verbatim}
\begin{equation*}
5
\end{equation*}
five prime numbers in the first two hundred
terms of the Smarandache prime sequence;
\begin{verbatim}
> nops(ssPrimes(primeSeq,200));
\end{verbatim}
\begin{equation*}
5
\end{equation*}
and two primes (11 and 1123) in the first one hundred and twenty
terms of the Smarandache Fibonacci sequence.
\begin{verbatim}
> ssPrimes(fibSeq,120);
\end{verbatim}
\begin{equation*}
[11, 1123]
\end{equation*}

It is clear that only the first term of the Smarandache even
sequence is prime. One interesting question, formulated in
\cite[Ch. 2]{MR99j:11005}, is the following:
\begin{quotation}
``How many elements of the Smarandache
even sequence are twice a prime?''
\end{quotation}
A simple search with \textsf{Maple}
shows that $2468101214$ is the only number twice a prime in the
first four hundred terms of the Smarandache even sequence (the
term 400 of the Smarandache even sequence is a number with 1147
decimal digits).
\begin{verbatim}
> ssTwicePrime := proc(n)
>   local ar, i:
>   ar := select(i->isprime(i/2),ss(evenSeq,n)):
>   convert(ar,list):
> end proc:
> ssTwicePrime(800);
\end{verbatim}
\begin{equation*}
[2468101214]
\end{equation*}


\section{Smarandache Relationships}

We now consider the so called \emph{Smarandache function}.
This function $S(n)$ is important for many reasons
(\textrm{cf.} \cite[pp. 91--92]{MR97m:01073}).
For example, it gives a necessary and sufficient condition for a number
to be prime: $p > 4$ is prime if, and only if, $S(p) = p$.
\emph{Smarandache numbers} are the values of the Smarandache function.

\subsection{Sequences of Smarandache numbers}

The \emph{Smarandache function} is defined in \cite{PropNumb} as follows:
$S(n)$ is the smallest positive integer number such that $S(n)!$ is divisible
by $n$.
This function can be defined in \textsf{Maple} by the following procedure:

\begin{verbatim}
> S:=proc(n)
>   local i, fact:
>   fact:=1:
>   for i from 2 while irem(fact, n)<>0 do
>     fact:=fact*i:
>   od:
>   return i-1:
> end proc:
\end{verbatim}
The first terms of the Smarandache sequence are easily obtained:
\begin{verbatim}
> seq(S(n),n=1..20);
\end{verbatim}
\begin{equation*}
 1, 2, 3, 4, 5, 3, 7, 4, 6, 5, 11, 4, 13, 7, 5, 6, 17, 6, 19, 5
\end{equation*}
\\

A sequence of $2k$ Smarandache numbers satisfy a Smarandache $k$-$k$
additive relationship if
\begin{equation*}
S(n)+S(n+1)+\dots+S(n+k-1)=S(n+k)+S(n+k+1)+\dots+S(n+2k-1) \, .
\end{equation*}
In a similar way, a sequence of $2k$ Smarandache numbers satisfy a
Smarandache $k$-$k$ subtractive relationship if
\begin{equation*}
S(n)-S(n+1)-\dots-S(n+k-1)=S(n+k)-S(n+k+1)-\dots-S(n+2k-1) \, .
\end{equation*}

In \cite{MR1764401,oai:arXiv.org:math/0010132} one finds the following
questions:
\begin{quotation}
``How many quadruplets verify a Smarandache $2$-$2$ additive relationship?''
\end{quotation}
\begin{quotation}
``How many quadruplets verify a Smarandache $2$-$2$ subtractive relationship?''
\end{quotation}
\begin{quotation}
``How many sextuplets verify a Smarandache $3$-$3$ additive relationship?''
\end{quotation}
To address these questions, we represent
each of the relationships by a \textsf{Maple} function:

\begin{verbatim}
> add2_2:=(V,n)->V[n]+V[n+1]=V[n+2]+V[n+3]:
> sub2_2:=(V,n)->V[n]-V[n+1]=V[n+2]-V[n+3]:
> add3_3:=(V,n)->V[n]+V[n+1]+V[n+2]=V[n+3]+V[n+4]+V[n+5]:
\end{verbatim}

We compute the first $10005$ numbers of the Smarandache sequence:

\begin{verbatim}
> SSN:=[seq(S(i),i=1..10005)]:
\end{verbatim}
\begin{verbatim}
(elapsed time: 59m29s)
\end{verbatim}

With the following procedure, we can identify all the positions in the sequence
$V$ that verify the relationship $F$.

\begin{verbatim}
> verifyRelation:=proc(F,V)
>   local i, VR: VR:=[]:
>   for i to nops(V)-5 do
>     if F(V,i) then VR:=[VR[], i]: fi:
>   od:
>   return VR;
> end proc:
\end{verbatim}
We can answer the above mentioned questions for the first
$10000$ numbers of the Smarandache sequence.

The positions verifying the Smarandache $2$-$2$ additive relationship are:
\begin{verbatim}
> V1:=verifyRelation(add2_2,SSN);
\end{verbatim}
\begin{equation*}
V1 := [6, 7, 28, 114, 1720, 3538, 4313, 8474]
\end{equation*}

Similarly, we determine the positions verifying the
Smarandache $2$-$2$ subtractive relationship,
\begin{verbatim}
> V2:=verifyRelation(sub2_2,SSN);
\end{verbatim}
\begin{equation*}
V2 := [1, 2, 40, 49, 107, 2315, 3913, 4157, 4170]
\end{equation*}
and the positions verifying the Smarandache $3$-$3$ additive relationship:
\begin{verbatim}
> V3:=verifyRelation(add3_3,SSN);
\end{verbatim}
\begin{equation*}
V3 := [5, 5182, 9855]
\end{equation*}

The quadruplets associated with the positions \texttt{V1}
($2$-$2$ additive relationship) are given by
\begin{verbatim}
> map(i->printf("S(%d)+S(%d)=S(%d)+S(%d) [%d+%d=%d+%d]\n",
  i,i+1,i+2,i+3,S(i),S(i+1),S(i+2),S(i+3)), V1):
\end{verbatim}
\begin{verbatim}
S(6)+S(7)=S(8)+S(9) [3+7=4+6]
S(7)+S(8)=S(9)+S(10) [7+4=6+5]
S(28)+S(29)=S(30)+S(31) [7+29=5+31]
S(114)+S(115)=S(116)+S(117) [19+23=29+13]
S(1720)+S(1721)=S(1722)+S(1723) [43+1721=41+1723]
S(3538)+S(3539)=S(3540)+S(3541) [61+3539=59+3541]
S(4313)+S(4314)=S(4315)+S(4316) [227+719=863+83]
S(8474)+S(8475)=S(8476)+S(8477) [223+113=163+173]
\end{verbatim}
We remark that in M.~Bencze's paper \cite{MR1764401}
only the first three quadruplets were found.
The quadruplets associated with the positions \texttt{V2}
($2$-$2$ subtractive relationship) are:
\begin{verbatim}
> map(i->printf("S(%d)-S(%d)=S(%d)-S(%d) [%d-%d=%d-%d]\n",
  i,i+1,i+2,i+3,S(i),S(i+1),S(i+2),S(i+3)), V2):
\end{verbatim}
\begin{verbatim}
S(1)-S(2)=S(3)-S(4) [1-2=3-4]
S(2)-S(3)=S(4)-S(5) [2-3=4-5]
S(40)-S(41)=S(42)-S(43) [5-41=7-43]
S(49)-S(50)=S(51)-S(52) [14-10=17-13]
S(107)-S(108)=S(109)-S(110) [107-9=109-11]
S(2315)-S(2316)=S(2317)-S(2318) [463-193=331-61]
S(3913)-S(3914)=S(3915)-S(3916) [43-103=29-89]
S(4157)-S(4158)=S(4159)-S(4160) [4157-11=4159-13]
S(4170)-S(4171)=S(4172)-S(4173) [139-97=149-107]
\end{verbatim}
Only the first two and fourth quadruplets were found
in \cite{MR1764401}.
The following three sextuplets verify a Smarandache
$3$-$3$ additive relationship:
\begin{verbatim}
> map(i->printf("S(%d)+S(%d)+S(%d)=S(%d)+S(%d)+S(%d)
  [%d+%d+%d=%d+%d+%d]\n",i,i+1,i+2,i+3,i+4,i+5,
  S(i),S(i+1),S(i+2),S(i+3),S(i+4),S(i+5)), V3):
\end{verbatim}
\begin{verbatim}
S(5)+S(6)+S(7)=S(8)+S(9)+S(10) [5+3+7=4+6+5]
S(5182)+S(5183)+S(5184)=S(5185)+S(5186)+S(5187) [2591+73+9=61+2593+19]
S(9855)+S(9856)+S(9857)=S(9858)+S(9859)+S(9860) [73+11+9857=53+9859+29]
\end{verbatim}
Only the first sextuplet was found by M.~Bencze's in \cite{MR1764401}.
For a deeper analysis of these type of relationships, see
\cite{subtRelation,MR1821195}.

\subsection{An example of a Smarandache partial perfect additive sequence}

Let $\{a_n\}$, $n\geq1$, be a sequence constructed in the following way:
\begin{equation*}
\begin{array}{l}
a_1=a_2=1;\\
a_{2p+1}=a_{p+1}-1;\\
a_{2p+2}=a_{p+1}+1 \, .
\end{array}
\end{equation*}
The following \textsf{Maple} procedure defines $a_n$.
\begin{verbatim}
> a:=proc(n)
>   option remember:
>   if (n=1) or (n=2) then return 1:
>   elif type(n, odd) then return a((n-1)/2+1)-1:
>   else return a((n-2)/2+1)+1:
>   fi:
> end proc:
\end{verbatim}
In \cite{MR1764401} the first 26 terms of the sequence are presented as being
\begin{verbatim}
> A:=1,1,0,2,-1,1,1,3,-2,0,0,2,1,1,3,5,-4,-2,-1,1,-1,1,1,3,0,2:
\end{verbatim}
One easily concludes, as mentioned in \cite{MR1821197},
that starting from the thirteenth term
the above values are erroneous. The correct values
are obtained with the help of our procedure:
\begin{verbatim}
> seq(a(i),i=1..26);
\end{verbatim}
\begin{equation*}
 1, 1, 0, 2, -1, 1, 1, 3, -2, 0, 0, 2, 0, 2, 2, 4, -3, -1, -1, 1, -1, 1, 1, 3,
 -1, 1
\end{equation*}
We prove, for $1 \leq p \leq 5000$, that
$\{a_n\}$ is a Smarandache partial perfect additive sequence,
that is, it satisfies the relation
\begin{equation}
\label{eq:sppasproperty}
a_1+a_2+ \dots +a_p = a_{p+1}+a_{p+2}+ \dots +a_{2p}  \, .
\end{equation}
This is accomplished by the following \textsf{Maple} code:
\begin{verbatim}
> sppasproperty:=proc(n)
>   local SPPAS, p;
>   SPPAS:=[seq(a(i),i=1..n)];
>   for p from 1 to iquo(n,2) do
>     if evalb(add(SPPAS[i], i=1..p)<>add(SPPAS[i], i=p+1..2*p))
>        then return false;
>     fi;
>   od;
>   return true;
> end proc:
\end{verbatim}
\begin{verbatim}
> sppasproperty(10000);
\end{verbatim}
\begin{equation*}
true
\end{equation*}
\begin{verbatim}
(elapsed time: 11.4 seconds)
\end{verbatim}
We remark that the erroneous sequence $A$
does not verify property \eqref{eq:sppasproperty}.
For example, with $p = 8$ one gets:
\begin{verbatim}
> add(A[i],i=1..8)<>add(A[i],i=9..16);
\end{verbatim}
\begin{equation*}
8 \neq 10
\end{equation*}


\section{Other Smarandache Definitions and Conjectures}

The Smarandache prime conjecture share resemblances
(a kind of dual assertion) with the famous Goldbach's conjecture:
``Every even integer greater than four can be expressed as
a sum of two primes''.

\subsection{Smarandache Prime Conjecture}

In \cite{MR1764401,oai:arXiv.org:math/0010132,PropNumb}
the so called \emph{Smarandache Prime Conjecture}
is formulated: ``Any odd number can be expressed as the sum of two primes
minus a third prime (not including the trivial solution $p=p+q-q$ when
the odd number is the prime itself)''.

We formulate a strong variant of this conjecture, requiring the odd number
and the third prime to be different (not including the situation $p=k+q-p$),
that is, we exclude the situation addressed by Goldbach's conjecture
(where the even integer $2p$ is expressed as the sum of two primes $k$ and $q$).

The number of times each odd number can be expressed as the sum
of two primes minus a third prime, are called
\emph{Smarandache prime conjecture numbers}.
It seems that none of them are known (\textrm{cf.} \cite{MR1764401}).
Here we introduce the notion of
\emph{strong Smarandache $n$-prime conjecture numbers}:
the number of possibilities that each positive odd number
can be expressed as a sum of two primes minus a third
prime, excluding the trivial solution and
imposing our requirement that the odd number
and the third prime must be different,
using all possible combinations of the first $n$ primes.

Given $n$, the next procedure determines such numbers
for all positive odd integers less or equal than $lim$.

\begin{verbatim}
> spcn:=proc(lim, n)
>   local y, z, i, primos, num, mat:
>   mat:=array(1..lim, 1..2,[seq([`?`, 0], i=1..lim)]):
>   primos:=seq(ithprime(i), i=1..n);
>   for i from 1 to n do
>     for y in [primos[i..n]] do
>       for z in [primos] do
>         num:=primos[i]+y-z;
>         if (num>=1 and num<=lim and type(num, odd) and
>             z<>primos[i] and z<>y and z<>num) then
>           if mat[num, 2]=0 then mat[num, 1]:=[primos[i], y, z]:
>           fi:
>           mat[num, 2]:=mat[num, 2]+1;
>         fi:
>       od:
>     od:
>   od:
>   for i by 2 to lim do
>     if mat[i, 2]=0 then printf("%d=? (0 possibilities)\n", i):
>     else printf("%d=%d+%d-%d (%d possibilities)\n", i,
>                 op(mat[i, 1]), mat[i, 2]):
>     fi:
>   od:
>   evalm(mat):
> end proc:
\end{verbatim}
All positive odd numbers less or equal than $19$ can be expressed according
to the conjecture, using only the first six primes.\footnote{For each number, only
one of the possibilities is showed.}
\begin{verbatim}
> spcn(19,6):
\end{verbatim}
\begin{verbatim}
1=2+2-3 (6 possibilities)
3=5+5-7 (3 possibilities)
5=3+13-11 (2 possibilities)
7=5+5-3 (2 possibilities)
9=3+11-5 (7 possibilities)
11=3+13-5 (3 possibilities)
13=5+11-3 (2 possibilities)
15=5+13-3 (5 possibilities)
17=7+13-3 (3 possibilities)
19=11+11-3 (3 possibilities)
\end{verbatim}
\begin{verbatim}
(elapsed time: 0.0 seconds)
\end{verbatim}
As expected, if one uses the first 100 primes, the number of distinct possibilities,
for which each number can be expressed as in our conjecture, increases.
\begin{verbatim}
> spcn(19,100):
\end{verbatim}
\begin{verbatim}
1=2+2-3 (1087 possibilities)
3=5+5-7 (737 possibilities)
5=3+13-11 (1015 possibilities)
7=3+17-13 (1041 possibilities)
9=3+11-5 (793 possibilities)
11=3+13-5 (1083 possibilities)
13=3+17-7 (1057 possibilities)
15=3+17-5 (770 possibilities)
17=3+19-5 (1116 possibilities)
19=3+23-7 (1078 possibilities)
(elapsed time: 1.8 seconds)
\end{verbatim}
How many odd numbers less or equal to $10000$
verify the conjecture?\footnote{In the follow
\emph{spcn} procedure calls, we removed from its definition the last \emph{for}
loop (\emph{spcn} without screen output).}
\begin{verbatim}
> SPCN1:=spcn(10000,600):
\end{verbatim}
\begin{verbatim}
(elapsed time: 30m59s)
\end{verbatim}
\begin{verbatim}
> n:=0: for i by 2 to 10000 do if SPCN1[i,2]>0 then n:=n+1; fi; od: n;
\end{verbatim}
\begin{equation*}
4406
\end{equation*}
Using the first $600$ primes, only $4406$ of the $5000$ odd numbers verify the
conjecture. And if one uses the first $700$ primes?
\begin{verbatim}
> SPCN2:=spcn(10000,700):
\end{verbatim}
\begin{verbatim}
(elapsed time: 49m34s)
\end{verbatim}
\begin{verbatim}
> n:=0:
> for i by 2 to 10000 do if SPCN2[i,2]>0 then n:=n+1; fi; od;
> n;
\end{verbatim}
\begin{equation*}
5000
\end{equation*}

Using the first $700$ primes, all the odd numbers up to 10000 verify the conjecture.
We refer the readers interested in the Smarandache prime conjecture
to \cite{MR1821203}.

\subsection{Smarandache Bad Numbers}
``There are infinitely many numbers that cannot be expressed as the difference
between a cube and a square (in absolute value). They are called
\emph{Smarandache Bad Numbers}(!)'' -- see \cite{MR1764401}.\\

The next procedure determines if a number $n$ can be expressed in the form
$n=|x^3-y^2|$ (\textrm{i.e.}, if it is a non Smarandache bad number),
for any integer $x$ less or equal than $x_{max}$.
The algorithm is based in the following equivalence
\begin{equation*}
n=|x^3-y^2| \quad \Leftrightarrow \quad y=\sqrt{x^3-n}\quad \vee \quad
y=\sqrt{x^3+n} \, .
\end{equation*}
For each $x$ between $1$ and $x_{max}$, we try to find an \emph{integer} $y$
satisfying $y=\sqrt{x^3-n}$ or $y=\sqrt{x^3+n}$,
to conclude that $n$ is a non Smarandache bad number.
\begin{verbatim}
> nsbn:=proc(n,xmax)
>   local x, x3:
>   for x to xmax do
>     x3:=x^3;
>     if issqr(x3-n) and x3<>n then return n[x, sqrt(x3-n)];
>     elif issqr(x3+n) then return n[x, sqrt(x3+n)]; fi;
>   od:
>   return n[`?`, `?`]
> end proc:
\end{verbatim}
F.~Smarandache \cite{PropNumb} conjectured that the numbers
$5,6,7,10,13,14,\dots$ are probably bad numbers.
We now ask for all the non Smarandache bad numbers which are less
or equal than $30$, using only the $x$ values between $1$ and
$19$. We use the notation $n_{x,y}$ to mean that $n=|x^3-y^2|$.
For example, $1_{2,3}$ means that $1=|2^3-3^2| = |8-9|$.

\begin{verbatim}
> NSBN:=map(nsbn,[$1..30],19);
\end{verbatim}
\begin{gather*}
NSBN:=[1_{2,3},2_{3,5},3_{1,2},4_{2,2},5_{?,?},6_{?,?},7_{2,1},8_{1,3},9_{3,6},
10_{?,?},\\
11_{3,4},12_{13,47},13_{17,70},14_{?,?},15_{1,4},16_{?,?},17_{2,5},18_{3,3},
19_{5,12},20_{6,14},\\
21_{?,?},22_{3,7},23_{3,2},24_{1,5},25_{5,10},26_{3,1},27_{?,?},28_{2,6},
29_{?,?},30_{19,83}]
\end{gather*}
As proved by Maohua Le in \cite{BadNumb}, we have just shown that
$7$ and $13$ are non Smarandache bad numbers:
$7=|2^3-1^2|$ and $13=|17^3-70^2|$.
The possible Smarandache bad numbers are:

\begin{verbatim}
> select(n->evalb(op(1,n)=`?`), NSBN);
\end{verbatim}
\begin{equation*}
[5_{?,?},6_{?,?},10_{?,?},14_{?,?},16_{?,?},21_{?,?},27_{?,?},29_{?,?}]
\end{equation*}

Finally, we will determine if any of these eight numbers
is a non Smarandache bad number, if one
uses all the $x$ values up to $10^8$.

\begin{verbatim}
> map(nsbn,[5,6,10,14,16,21,27,29],10^8);
\end{verbatim}
\begin{gather*}
[5_{?,?},6_{?,?},10_{?,?},14_{?,?},16_{?,?},21_{?,?},27_{?,?},29_{?,?}]
\end{gather*}
\begin{verbatim}
(elapsed time: 14h30m)
\end{verbatim}
From the obtained result, we conjecture that
$5,6,10,14,16,21,27$, and $29$, are bad numbers.
We look forward to readers explorations and discoveries.


\section*{Acknowledgments}

This work was partially supported by the
program PRODEP III/5.3/2003, and by the
R\&D unit CEOC of the University of Aveiro.



\end{document}